\documentclass[12pt]{amsart}
\usepackage{amssymb}
\usepackage[all]{xy}
\newdir{ >}{{}*!/-5pt/@{>}}
\newtheorem{theorem}{Theorem}
\newtheorem{proposition}[theorem]{Proposition}
\newtheorem{corollary}[theorem]{Corollary}
\newtheorem{lemma}[theorem]{Lemma}
\theoremstyle{remark}
\newtheorem{remark}[theorem]{Remark}
\newtheorem{example}[theorem]{Example}
\newtheorem{definition}[theorem]{Definition}
\title[Mapping spaces between manifolds]{Hochschild cohomology of Sullivan algebras and mapping spaces between manifolds}
\author{J.-B. Gatsinzi}
\address{Department of Mathematics and Statiscal Sciences, Botswana International University of Science and Technology.}
\email{gatsinzij@biust.ac.bw}
\thanks{A partial support from the IMU-Simons Africa Fellowship is acknowledged}

\DeclareMathOperator{\Ext}{Ext}
\DeclareMathOperator{\Hom}{Hom}
\DeclareMathOperator{\Der}{Der}
\DeclareMathOperator{\map}{map}

\DeclareMathOperator{\aut}{aut}

\begin{document}

	\begin{abstract}
		Let $e: N^n \rightarrow M^m $ be an embedding into a compact manifold $M$. We study the relationship between the   homology of the free loop space $LM$ on $M$  and  of the space $L_NM$ of loops of $M$ based in $N$ and define a shriek   map
		  $ e_{!}: H_*( LM, \mathbb{Q}) \rightarrow   H_*( L_NM, \mathbb{Q})$ using Hochschild cohomology and study its properties. We  also extend a result   of F\'elix on the injectivity of the induced map $ \aut_1M  \rightarrow \map(N, M; f ) $  on rational homotopy groups when  $M$ and $N$ have the same dimension and $ f: N\rightarrow M $ is a map of  non zero degree.   
	\end{abstract}
	\maketitle 
	\section{Introduction}
 
 All spaces are assumed to be simply connected and  (co)homology  coefficients are taken in the field $ \mathbb{Q} $ of rationals. If $M$ is a compact oriented manifold of dimension $m$ and $ LM= \map(S^1, M)$ the space of free loops  in $M$, then there is an intersection product 
 $$ \mu: H_{p+m} (LM) \otimes H_{q+m} (LM)\rightarrow H_{p+q+m} (LM) $$ which induces  a graded  multiplication on $ \mathbb{H}_*(LM) = H_{*+m} (LM)$ turning into a graded algebra~\cite{c-s}. Consider  the embedding $ e: N \rightarrow M$ of  a submanifold of degree $n$. Construct the pullback
 $$
 \xymatrix{L_NM \ar[r]^{\tilde e} \ar[d]_{\tilde{p}} & LM  \ar[d]^{p}\\
 	N \ar[r]_{e}  & M,}
$$
where $p$ is the evaluation of a loop at $1 \in S^1$.
There is also an intersection product   on $\mathbb{H}_*(L_NM) = H_{*+n} (L_NM)$, turning it into  commutative graded algebra~\cite{sul-2004}. \\

We  consider a morphism $ f: (A,d) \rightarrow (B,d) $  of commutative differential graded algebras which models the embedding $e$, where $(A,d) $ and $ (B,d) $ are Poincar\'e duality algebras~\cite{felix-2011}.   We show that there is an  $A$-linear shriek map $f_{!}: (B,d) \rightarrow (A,d) $   of degree $m-n$.
We consider induced maps  $ HH^*(f) : HH^*(A ; A)  \rightarrow HH^*(A;  B)$  and $HH^*(f_{!}): HH^*(A;B ) \rightarrow HH^*(A; A) $ in Hochschild cohomology. Moreover we obtain the following.
\begin{theorem}
	\label{first}
	The composition  map $$  HH^*(f_{!}) \circ HH^*(f):  HH^*(A; A)  \rightarrow HH^*(A; A)$$  is the multiplication by the Poincar\'e dual of the fundamental class of $N$ in $M$. 
\end{theorem}

\begin{theorem}
	\label{second}
Let $g: N^m \rightarrow M^m $ be a map between manifolds of same dimension $m$ such $ \deg f \neq 0$
and $f: (A,d ) \rightarrow (B,d) $ a cdga model of $g$. Then
$$ HH^*(A; A) \rightarrow HH^*(A; B )$$
is injective. 
\end{theorem}
The above result suggests that $\mathbb{H}(\tilde g): \mathbb{H}_*(L_NM) \rightarrow \mathbb{H}_*(LM) $ is an injective algebra homomorphism, where $ \tilde g: L_NM \rightarrow LM$ is the pullback of $g: N \rightarrow M$ along the fibration  $p: LM \rightarrow M$  defined by $ p(\gamma) = \gamma(0).$ \\

The paper is organized as follows: In Section~$2$ we define a shriek map  $f_{!}: (B,d) \rightarrow (A,d) $ and prove Theorem~\ref{first}. In Section~$3$, we recall a resolution to compute $ HH^*(C^*(M), C^*(N)) $ and in Section~$4$ we prove Theorem~\ref{second}.

\section{A shriek map}
We first recall some facts in Rational Homotopy Theory. We make use of Sullivan models for which the standard reference is \cite{fht-01}. All vector spaces are over the ground field  $\mathbb{Q}$.  A differential graded algebra $(A,d)$ is a direct sum of vector spaces $A^p$, that is, 
$ A = \oplus_{p\geq 0} A^p$  together with a graded multiplication
$ \mu : A^p \otimes A^q \rightarrow A^{p+q} $ which is associative.  An element $ a \in A^p$ is called homogeneous of degree $ |a| = p$. Moreover there is a differential $ d: A^p \rightarrow A^{p+1} $ which an algebra derivation, that is, $d(ab) = (da)b + (-1)^{|a|} a (db)$ and satisfies $ d^2 =0$. \\

The algebra $A$ is commutative if $ ab = (-1)^{|a||b|} ba $.  
If $(A,d) $ is a commutative differential graded algebra (cdga for short), then $ H^*(A,d) $ is graded commutative.   A morphism $ f: (A,d) \rightarrow (B,d) $ of cdga's  is called a quasi-isomorphism if $ H^*(f) $ is an isomorphism. 
A cdga $(A,d)$ is called  simply connected if  $ H^0(A) = \mathbb{Q}$ and  $ H^1(A) = 0$.\\

A commutative graded algebra $A$ is free if it is of the form $ \land V = S(V^{even})  \otimes E(V^{odd}) $, where $V^{even} = \oplus_{i \geq 1} V^{2i}$  and $ V^{odd} = \oplus_{i \geq 0} V^{2i+1}$. 
 A Sullivan algebra is a cdga $( \land V, d)$, where $V= \oplus_{i \geq 1} V^i $ admits a homogeneous basis $ \{ x_i\}_{ i\in I} $ indexed by a well ordered set $I$ such $ dx_i \in \land (\{x_i\})_{ i < j}$. A Sullivan algebra is called minimal if  $ dV \subset \land^{\geq 2}V$~\cite{fht-01}. If there is a quasi-isomorphism  $f:  (\land V, d)\rightarrow (A,d)$, where $ (\land V, d)$ is a (minimal) Sullivan algebra, then we say that $ (\land V,d) $ is a (minimal) Sullivan model of $(A,d)$.
\\

 To a simply connected  topological space $X$ of finite type,  Sullivan associates in a functorial way a cdga $ A_{PL}(X) $ of piecewise linear forms on $X$ such  $H^*(A_{PL}(X)) \cong H^*(X, \mathbb{Q})$~\cite{sul77}.  A Sullivan model of $X$ is a Sullivan model of $ A_{PL}(X) $. Moreover any cdga $(A,d) $  is called a cdga-model of $X$ if there is  a sequence of quasi-isomorphisms
 $$ (A,d)\rightarrow (A_1, d) \leftarrow \ldots \rightarrow (A_{n-1}, d) \leftarrow  A_{PL}(X). $$

We state here the fundamental result of Sullivan algebras.
\begin{proposition}
If $ (A,d)	$ is a simply connected cdga   then  there is a minimal Sullivan algebra $ (\land V, d)$ together with a quasi-isomorphism  $ (\land V,d )\rightarrow (A,d) $. Moreover  $ (\land V, d)$  is unique up to isomorphism. It is called the   minimal Sullivan model of $(A,d)$~\cite[\S \ 12]{fht-01}.
\end{proposition}
\begin{definition}
	Let $X$ be a simply connected space. A  minimal Sullivan  model   $ (\land V, d) $ of $X$ is the   minimal Sullivan   model of $A_{PL}(X)$. It is called formal if there is a quasi-isomorphism $( \land V, d)\rightarrow H^*(\land V,d)$. In this case $X$ is called a formal space. Formal spaces include spheres, compact Lie groups and  complex projective spaces.
\end{definition}
\begin{definition} A commutative differential graded  algebra  $(A,d )$ is a   Poincar\'e algebra of formal dimension $n$ if   $A$ is connected and there is  a linear map $  \epsilon : A^n \rightarrow \mathbb{Q}$ such that 
\begin{enumerate}
	\item $ \epsilon (dA^{n-1}) = 0 $,
	\item the bilinear form $ b: A^k \otimes A^{n-k}\rightarrow \mathbb{Q} $, defined by $ b(x \otimes y) = \epsilon(ab)$ is non degenerate. 
\end{enumerate}

\end{definition}
\begin{remark}
	\label{rmk1}
	If $A$ is of finite type, then $ A^i = 0$ for $ i > n$ and $A$ is finite dimensional. Moreover if $ \{ a_1, \dotsc , a_k  \}$ is a homogeneous basis  of $A$, then there is a dual  homogeneous basis $ \{ a_j^*\} $ such that 
	$ \epsilon (a_i a_j^*) = \delta_{ij}$. We denote by $ a^{\#} $ the dual of $a$ in $ A^{\#} = \Hom (A, \mathbb{Q})$.  In particular $ \omega_A = \epsilon^{ \#} \in (A^{\#})^{\#} \cong A $ is the fundamental class of $A$. Moreover there is an isomorphism  of $A$-modules $ \pi_A: A \rightarrow A^{\#} $ defined by $\pi_A (a) (x)  = b(ax) .$
\end{remark}

If $ (\land V, d) $ is the minimal Sullivan model of a simply connected space $X$, where $ H^*(X, \mathbb{Q}) $ satisfies Poincar\'e duality, then $ (\land V, d)$  is quasi-isomorphic to a Poincar\'e duality algebra $(A,d) $~\cite{lamb-stanley}.
In particular, a simply connected smooth manifold $M$ of  dimension $m$ has a cdga-model $ (A,d)$ which satisfies Poincar\'e duality in dimension $m$. \\

Let $f: (A,d)\rightarrow (B,d)$ be a map between cdga's with Poincar\'e duality in dimensions $m$ and $n$ respectively.
We can now relate isomorphisms   $\pi_A:  A \stackrel{\simeq}{\rightarrow} A^{\#}$ and $ \pi_B : B \stackrel{\simeq}{\rightarrow}  B^{\#}$.
  
\begin{proposition} 
If $f$ is surjective, then 	there exists a morphism  of $A$-modules $ f_{!}: B \rightarrow A $  making the following diagram commutative.
$$
\xymatrix{
B  \ar[d]_{\simeq }^{\pi_B} \ar[r]^{f_{!}} &  A  \ar[d]^{\cong}_{\pi_A}\\
B^{\#}  \ar[r]^{f^{\#}}&   A^{\#}
}
$$
\label{prop}
\end{proposition}
\proof 

Let $1 \in B$, then $ \pi_B(1) = \omega_B^{\#}$, where $ \omega_B$ is a cocycle which represents the fundamental class $[ \omega_B] \in H^n(B)$. As $ \pi_A $ is bijective, there exists $ \alpha \in A$ such that $ \pi_A(\alpha) = f^{\#} (\omega_B^{\#})$. As  $f$ is surjective, then given $b \in B$, there exists $ a\in A$ such that $ b= f(a)$. Recall that $B$ is an $A$-module through the action induced by $f$, hence $ b = f(a)1 = a* 1$. Therefore we define $ f_{!} (b) = a\alpha$. In particular $ f_! f(a) =a \alpha $.  \\
We show that the above diagram commutes. Let $ b \in B$ and $ a\in A$ such that $ b= f(a)$.
On one hand
\begin{equation}
f^{\#} (\pi_B(b))  = f^{\#} (\pi_B(b \times 1))= f^{\#} (b \omega_B^{\#}),                  
                  \end{equation}
whereas
\begin{equation}
\pi_A(f_{!} (b))   = \pi_A (a \alpha)  
                =a \pi_A(\alpha) 
                = af^{\#}(\omega_B^{\#})  .
\end{equation}
Let $x\in A$. Then
\begin{equation} f^{\#}( b \omega_B^{\#}) (x)= (b \omega_B^{\#})(f(x)) = \omega_B^{\#} (b f(x)),
	\end{equation}
and
\begin{equation} 
\begin{array}{ll} (a f^{\# }( \omega_B^{\#})) (x)   & = (f^{\#}(\omega_B^{\#}))(ax) = \omega_B^{\#} (f(ax))  \\
 &   = \omega_B^{\#} (f(a)f(x)) = \omega_B^{\#} (bf(x)).
 \end{array}
\end{equation}
Hence  $ f^{\#} (b \omega_B^{\#}) =  af^{\#}(\omega_B) $ and the diagram commutes. \\

Finally we show that $ f_{!}$ is a morphism of $A$-modules. If $ x\in A$ and $ b\in B$, then 
$$ f_{!}(x* b)= f_{!} (f(x)b)= f_{!}(f(xa))= (xa) \alpha = x f_{!}(b).$$  In particular $ f_{!}(b) = f_{!} (b \times 1) = a* f_{!} (1)$.  \qed 

\begin{remark}
	\label{rmk2}
	If $ \omega_B$ is a cocycle representing the fundamental class of $(B,d)$ and $ f$ is surjective, then there exists $ x \in A$ such that $ f(x) =\omega_B$. Then $ f^{\#} (\omega_B^{\#}) = x^{\#} =\pi_A(x^*)$,  where $ x^{*}$  is the dual of $x$ under  a choice of a  basis  $ \{ a_i\}$ of $A$ and its dual $ \{ a^*_j\}  $ (see Remark~\ref{rmk1}). If $ dx = 0$, then $ [x] \in H^*(A) \neq 0$ and  $ [x^*]  \in H^{m-n} (A)$ is non zero.
\end{remark}

\begin{example}
	Consider the inclusion $ i: \mathbb{C}P^n \rightarrow \mathbb{C}P^{n+k} $. As complex projective spaces are formal, a  cdga model  of the inclusion is
	$$f: \land x_2/(x_2^{n+k+1})\rightarrow \land y_2/(y_2^{n+1}),  $$
	 where $ f(x) = y$. 
	Then $ f_{!} $ is defined by $ f_{!}(1)= x^k$. Hence $f_{!}(y^i) =  x^{k+i}, $ for $ 0 \leq i \leq n $.
	
\end{example}
 \section{Hochschild cohomology}
 If $(A,d) $ is a graded differential algebra and $ (M,d) $ a graded $A$-bimodule, then the Hochschild cohomology of $A$ with coefficients in $M$ is defined by $ HH^*(A; M) = \Ext_{A^e}(A, M)$, where $ A^e = A \otimes A^{opp}$.\\
 
 Let  $A= (\land  V, d)$ be the  minimal Sullivan model of a simply connected space $X$. Then   
 \begin{equation}
 \label{semifree-resolution}
  P=(\land V \otimes \land V \otimes \land \bar V, \tilde D)\rightarrow ( \land V, d)
  \end{equation}
  is a semi-free resolution  of $ \land V$ as a $\land V \otimes \land V $-module, where $ \bar V = sV$~\cite{felix-2011}.
   
 Moreover, the pushout
 $$
 \xymatrix{
 	(\land V \otimes \land V, d\otimes 1 + d \otimes 1) \ar[d]^{\mu}\ar@{>->}[r] & (\land V \otimes \land V \otimes \land \bar V, \tilde D) \ar[d] \\
 	(\land V, d)\ar@{>->}[r]   &   (\land V \otimes \land \bar V, D) 
 }
 $$
 yields a Sullivan model
 $ (\land V \otimes \land \bar V, D)$  of the free loop space on   $X$~\cite{s-v}.  The differential is given by $Dv = dv $ for $v \in V$ and $D \bar v  = -Sdv$, where $S$ is the unique derivation on $ \land V \otimes \land \bar V$ defined by $ Sv = \bar v $ and $ S \bar v = 0$. \\[1ex]
 Hence if  $(M,d)$ is a $ \land V$-differential module, then the Hochschild cochains $ CH(A;M) $ are given by
 \begin{equation}
 \label{HH-cochain}
 \begin{array}{ll}
 CH^*(A; M) & = (\Hom_{\land V \otimes \land V} ( \land V \otimes \land V \otimes \land \bar V, M), D)  \\
 & \cong (\Hom_{\land V} ( \land V \otimes \land \bar V, M), D).
 \end{array}
 \end{equation}
 As the differential of $ D$ on $ \land V \otimes \land \bar V$ satisfies $$ D (\land V \otimes \land^n \bar V)  \subset \land V \otimes \land^n \bar V,$$ one gets  a Hodge type decomposition
 $$HH^*(A;  M) = \oplus_{i \geq 0} HH_{(i)}^* (A; M), $$ where $  HH_{(i)}^* (A; M)= H^*(\Hom_{\land V} ( \land V\otimes  \land ^i \bar V, \land V), D)$. 
 Moreover, if $ L = s^{-1} \Der \land V $, then the symmetric  algebra  $  (\land _A L, d) $ 
 is quasi-isomorphic to the Hochschild cochain complex $ (\Hom_{\land V} ( \land V \otimes \land \bar V, \land V), D)$~\cite{gat10}.  
 If $ (\land V, d)$ the minimal Sullivan model of a simply connected smooth compact and oriented  manifold $M$ of dimension $m$, then there is an isomorphism of BV-algebras
 $\mathbb{H}_* (LM) \cong HH^*(\land V;  \land V)$ ~\cite{c-j, FTV-2004, ftv-08}. \\
 
 Let $M$ be a smooth  compact, oriented and simply connected manifold of dimension $m$.  For submanifolds $N$ and $N'$, we  denote by $L_N^{N'}M $ the space of paths in $M$ starting in $N$ and ending in $N'$, and $ L_N^N M$ is simply written $L_NM$. Let  $N_1, N_2$ and $N_3$ be submanifolds of $M$. 
  When coefficients are rationals (or in  $ \mathbb{Z}/n \mathbb{Z}$)  Sullivan showed that there is an intersection product
 $$ \mu: H_{p+d}(L_{N_1}^{N_2} M) \otimes H_{q+d} (L_{N_2}^{N_3} M)  \rightarrow  H_{p+q+d}(L_{N_1}^{N_3 } M)$$
  where $ d= \dim N_2$~\cite{sul-2004}. In particular if $ N_1=N_2=N_3 = N$, one gets a graded commutative algebra  structure on $ \mathbb{H}_*(L_N M , \mathbb{Q}) = H_{*+d} (L_N M, \mathbb{Q})$. \\

 Let $e: N^n \hookrightarrow M^m $ be an embedding where $N$ is simply connected  and $ f: (A,d) \rightarrow (B,d)$ a cdga model of $e$, where both $ (A,d) $ and $ (B,d)$ satisfy Poincar\'e duality.  Assume that $f$ is surjective and  let $[y] \in H^n(B) $  be the fundamental class. Let $  x \in A$ such that $f(x)=y$.  We will assume that $x$ is a cocycle and consider $ [x]  \in H^n(A,d) $. 
 
 \begin{theorem}
 	\label{thm1}
 Under the above hypotheses, the composition 
 $$ 
 \xymatrix{HH^*(A;A) \ar[r]^{HH^*(f)} & HH^*(A; B) \ar[r]^{HH^*(f_!)} & HH^*(A; A) }
 $$
 is the multiplication with the Poincar\'e dual $  [x^*] \in H^{m-n} (A,d) $  of $[x]$.
 \end{theorem}
 \proof We consider a minimal Sullivan model $ \phi: (\land V, d) \rightarrow (A,d)$. By Eq.~\eqref{HH-cochain}, $HH^*(A; A)$ is obtained as the cohomology of the complex 
 $$
 \begin{array}{ll}
  \Hom_{\land  V \otimes \land V} ( \land V \otimes \land V \otimes \land \bar V, \land V)   &\cong
  \Hom_{ \land V} ( \land V \otimes  \land \bar V, \land V) 
  \\
  &  \simeq \Hom_{ \land V} ( \land V \otimes   \land \bar V, A).
  \end{array}
  $$
  
  If $ \gamma \ \in \Hom_{ \land V} ( \land V \otimes   \land \bar V, A)$,
  then
$$ (CH(f_!)  \circ CH(f))( \gamma) (x)= (f_! \circ  f)(\gamma) (x) =   \alpha \gamma (x),
$$ 
where $ \alpha = x^*$, by Remark~\ref{rmk2}.
Therefore, if $ \gamma $ is a cocycle, then
$ HH^*(f_!) \circ HH^*(f) = [x^*][ \gamma]$. 
 \qed 
 
 \begin{example}
 The hypotheses of Theorem~\ref{thm1} are satisfied if $ e: N \rightarrow M$ is an embedding between formal smooth manifolds where $ H^*(e)$ is surjective, for instance the inclusion $ \mathbb{C}P^n\rightarrow \mathbb{C}P^{n+k}$.  Let $ A= H^*(\mathbb{C}P^{n+k}, \mathbb{Q}) = \land x_2/(x_2^{n+k+1})$. The complex to compute $HH^*(A;A) $ is given by
 	$ (A \otimes \land(z_1, z_{2(n+k)}), D) $ where subscripts indicate the lower degree, and $ Dz_{2(n+k)} =0$ and $ Dz_1= (n+k+1) x_2^{n+k}z_{2(n+k)}$~\cite{gat16}. Here an element $x\in A^n=A_{-n}$ is assumed to be of lower degree $-n$.  At  chain's level, the composition
 	$$ CH^*(f_!)\circ CH(f):  (A \otimes \land(z_1, z_{2(n+k)}), D) \rightarrow  (A \otimes \land(z_1, z_{2(n+k)}), D)$$
 	is the multiplication by $ x^k$.
 \end{example}
 If $e:  N \rightarrow M$ is an embedding between manifolds,
then $L_NM $ 
is the pullback of the following diagram
\begin{equation}
\label{loop1}
\xymatrix{
	L_NM  \ar[r]^{\tilde e}  \ar[d]^{\tilde p}&    LM \ar[d]^{p}\\
	N \ar[r]^{e}  &   M,
}
\end{equation}
where $ p(\gamma) = \gamma (0)$.\\

Assume that $\pi_*(M) \otimes \mathbb{Q} $ is finite dimensional and $ (\land V, d) $ is the minimal Sullivan model of $M$. Then $HH^*(\land V; \land V)$ is the homology of the  complex $ (\land V \otimes  \land Z, D)$ where $ Z \simeq s^{-1} V^{\#}$~\cite{gat16}.
\begin{proposition}
  If $ f: (A,d) \rightarrow (B,d) $ is a model of $e: N \rightarrow M$, then $HH^*(C^*(M); C^*(N))$ is computed by the complex  $ (B \otimes \land Z, D)   $ obtained as the pushout
\begin{equation}
\xymatrix{
(A,d) \ar[r] \ar[d]  &    (A \otimes \land Z, D) \ar[d] \\
(B,d) \ar[r] &    (B \otimes \land Z, D)
}
\end{equation}
\proof  
Let  $ (\land V, d)$ is the minimal Sullivan model of $M$, where $V$ is finite dimensional. Then $ \mathbb{H}_*(LM) $ is the homology of the complex  $ (\land V \otimes \land Z, D)$ where $ Z = s^{-1} V^{\#}$ and the differential $D$ is induced by $ \delta $ on $ (\Der \land V, \delta)$ where  $ V^{\#} \subset \Der \land V$.  
As $ (\land V, D) \rightarrow (A,d) $ is a quasi-isomorphism, then the pushout
is a model of the pullback in Eq.~\ref{loop1}. \qed 

However, it is known whether  structure of $\mathbb{H}_*(L_NM) $ and $  H_*(B \otimes  \land Z, D)$ are isomorphic as algebras.
\end{proposition} 
 \section{Maps between manifolds of same dimension}
 Let $ f: (A,d) \rightarrow (B,d)$ be a morphism  of graded cochain algebras.  An $f$-derivation of degree $n$ is a linear map $\theta:  A^* \rightarrow B^{*-n}$ such that $ \theta (xy) = \theta(x)f (y) + (-1)^{n|x|}f (x) \theta (y)$.
 We denote by $ \Der_n(A, B; f) $ the vector space of all $ f$-derivations of degree $n$ and $ \Der(A, B ; f ) = \oplus_n \Der_n(A, B; f)$. 
 Define a differential $ \delta  $ on $ \Der(A, B ; f )$ by $ \delta \theta = d_B \theta - (-1)^{|\theta|} \theta d_A$.
 If $A= B$, then we simply write $ \Der A$ for $  \Der(A, A ; 1_A ) $.  The graded vector space $ \Der A$ is endowed with the commutator bracket  turning it into a graded differential Lie algebra.  There is an action of $A$ on $ \Der A$, defined by $ (a\theta)(x) = a \theta (x)$, making $ (\Der A, \delta ) $ a differential graded module over $A$.  \\

 Let $M$ and $N$ be compact and oriented  manifolds of dimension $n$ and $ g: N \rightarrow M $  a smooth map such that $ \deg g  \neq 0$.  Consider a Poincar\'e duality  model $ f : (A,d)\rightarrow (B,d)$   of $g$.
 Then $f$ is injective and $ B = f(A) \oplus Z$, where $ dZ \subseteq Z $ and $ f(A).Z$~\cite{felix-2011}. Therefore $Z$ is an $A$-submodule. Moreover the projection 
 $p: B = f(A) \oplus Z \rightarrow A$ is a morphism of $A$-modules.  
 \begin{theorem}[\cite{felix-2011}, Theorem 2]
 	 \label{injection}
 Consider a  surjective  Sullivan model $ \phi: (\land V, D ) \rightarrow (A,d)$. 
Then 
 \begin{equation}f_*: (\Der (\land V, A; \phi), \delta) \rightarrow ( \Der (\land V, B ;  f\circ \phi ), \delta) 
\end{equation}
  induces an injective map in homology.
  \end{theorem}
 This can be interpreted in terms of rational homotopy groups of function spaces. Let $ g: X \rightarrow Y$ be a continuous map  between CW complexes where $Y$ is finite and $X$ of finite type and $ \phi: (\land Z, d)  \rightarrow (B,d)$ a Sullivan model of $g$.  Consider $ \map(X, Y; g)$  be the space of continuous mappings from $X$ to $Y$ which are homotopic to $f$. 
 There is a natural isomorphism~\cite{block-lazarev, B-M2008, lupton-s-2007b}
 $$ \pi_n(\map(X, Y; g)) \otimes \mathbb{Q} \cong H_n(\Der(\land V, B ; \phi ), \delta ), \ n \geq 2.$$

 Hence if $ g: N \rightarrow M $ is a map between simply connected smooth manifolds such that $ \deg g \neq 0$, then  
 the  map
 $$j_M: \aut_1 M=\map(M, M; 1_M) \rightarrow \map(N, M; g)$$
 induces an injective map 
 $$ \pi_*(j_M) \otimes \mathbb{Q}: \pi_*(\aut_1 M) \otimes \mathbb{Q}\rightarrow \pi_*(\map(N,M;g))\otimes \mathbb{Q}. $$
Let $ \phi: (\land V, d) \rightarrow (A,d)$  be a Sullivan model and $ \rho = f \circ \phi$.
We have the following  commutative diagram
$$
\xymatrix{ H_* (\Der \land V, \delta)  \ar@{>->}[d] \ar@{^{(}->}[r] & HH^*(A; A) \ar[d] \\
	H_*(\Der (\land V, B ;\rho), \delta) \ar@{^{(}->}[r] &   HH^*(A; B),
}
$$ 
where
horizontal maps are inclusions~\cite{gat-2019}. 
We  show that the remaining vertical arrow is injective, which is a generalization of Theorem~\ref{injection}.
 
 \begin{theorem}
 	\label{thm-mfd}
 	Let $g: N \rightarrow M$ be a smooth map between manifolds and $ f: (A,d) \rightarrow (B,d) $ a Poincar\'e duality model of $g$. Then
 	the induced map
 	$$ 
 	\xymatrix{ HH^*(A; A)  \ar[r]^{HH^*(f)}& HH^*(A; B)
 	}
 $$
 is injective.
 \end{theorem}

 \proof 
   As $ B = f(A) \oplus Z$, then $f(A)= \rho (\land V) $ is a submodule of $B$ viewed as a $ \land V$-module and  $Z$ is also a  $ \land V$-submodule of $B$.  Therefore
 $$ \Hom_{\land V }(\land V \otimes \land \bar V, B) \cong  \Hom_{\land V }(\land V \otimes \land \bar V, f (A)) \oplus \Hom_{\land V }(\land V \otimes \land  \bar V, Z).$$
 Moreover, the projection $ p: B = f (A) \oplus Z \rightarrow f(A) \cong A $ is a morphism of $ \land V$-modules.  It induces a chain map
 $$ p_*: \Hom_{\land V }(\land V \otimes \land \bar V, B) \rightarrow  \Hom_{\land V }(\land V \otimes \land \bar V, A) $$
 such that $ p_* \circ f_* $ is the identity. Therefore $ f_* $ is injective in homology.  \qed

  We can then deduce the following
 \begin{corollary}
Under the hypotheses of Theorem~\ref{thm-mfd}, there is an injective  map $ H_*(f)_{!}: H_*(LM, \mathbb{Q}) \rightarrow H_* (L_NM, \mathbb{Q})  $
 \end{corollary}
\proof
Recall that  there is an isomorphism $ HH_*(A, A)  \cong H^*(LM)$~\cite{jones}. Dualizing this isomorphism  and using Poincar\'e duality yield 
an isomorphism $ HH^*(A; A^{\#}) \cong H_*(LM)$. In the same way, there is an isomorphism $ HH^*(A, B^{\#}) \cong H_*(L_NM)$, Then $H_*(f)_{!}$ is given by the composition
$$
\xymatrix{ HH^*(A; A^{\#})  \ar[r]^{(\pi_A)_*^{-1}} & HH^*(A;A)\ar[r]^{f_*} & HH^*(A;B) \ar[r]^{(\pi_B)_*} & HH^*(A; B^{\#}).
}
$$
Hence 	it is injective. \qed
\bibliographystyle{amsplain}
\bibliography{references} 

\providecommand{\bysame}{\leavevmode\hbox to3em{\hrulefill}\thinspace}
\providecommand{\MR}{\relax\ifhmode\unskip\space\fi MR }
\providecommand{\MRhref}[2]{%
  \href{http://www.ams.org/mathscinet-getitem?mr=#1}{#2}
}
\providecommand{\href}[2]{#2}
\begin{thebibliography}{10}

\bibitem{block-lazarev}
J.~Block and A.~Lazarev, \emph{Andr\'{e}-{Q}uillen cohomology and rational
  homotopy of function spaces}, Adv. Math. \textbf{193} (2005), 18--39.

\bibitem{B-M2008}
U.~Buijs and A.~Murillo, \emph{The rational homotopy {L}ie algebra of function
  spaces}, Comment. Math. Helv. \textbf{83} (2008), 723--739.

\bibitem{c-s}
M.~Chas and D.~Sullivan, \emph{String topology}, preprint math GT/9911159,
  1999.

\bibitem{c-j}
R.L Cohen and J.D.S Jones, \emph{A homotopy theoretic realisation of string
  topology}, Math. Ann. \textbf{324} (2002), no.~4, 773--798.

\bibitem{felix-2011}
Y.~F\'elix, \emph{Mapping spaces between manifolds and the evaluation map},
  Proc. Amer. Math. Soc. \textbf{139} (2011), 3763--3768.

\bibitem{fht-01}
Y.~F{\'e}lix, S.~Halperin, and J.-C. Thomas, \emph{Rational {H}omotopy
  {T}heory}, Graduate Texts in Mathematics, no. 205, Springer-Verlag, New-York,
  2001.

\bibitem{ftv-08}
Y.~F{\'e}lix, J.-C. Thomas, and M.~Vigu{\'e}, \emph{Rational string topology},
  J. Eur. Math. Soc. (JEMS) \textbf{9} (2008), 123--156.

\bibitem{FTV-2004}
Y.~F\'elix, J-C. Thomas, and M.~Vigu\'e-Poirrier, \emph{The {H}ochschild
  cohomology of a closed manifold}, Publ. Math. Inst. Hautes {\'E}tudes Sci.
  \textbf{99} (2004), 235--252.

\bibitem{gat10}
J.-B. Gatsinzi, \emph{Derivations, {H}ochschild cohomology and the {G}ottlieb
  group}, Homotopy Theory of Function Spaces and Related Topics (Y.~F{\'e}lix,
  G.~Lupton, and S.~Smith, eds.), Contemporary Mathematics, vol. 519, American
  Mathematical Society, Providence, 2010, pp.~93--104.

\bibitem{gat16}
J.-B Gatsinzi, \emph{Hochschild cohomology of a {S}ullivan algebra}, Mediterr.
  J. Math. \textbf{13} (2016), 3765--3776.

\bibitem{gat-2019}
J.-B. Gatsinzi, \emph{Hochschild cohomology of {S}ullivan algebras and mapping
  spaces}, Arab J. Math. Sci. \textbf{25} (2019), 123--129.

\bibitem{jones}
J.~D.~S. Jones, \emph{Cyclic homology and equivariant homology}, Inv. Math.
  \textbf{87} (1987), 403--423.

\bibitem{lamb-stanley}
P.~Lambrechts and D.~Stanley, \emph{Poincar\'e duality and commutative
  differential graded algebras}, Ann. Sci. {\'E}c. Norm. Sup{\'e}r. \textbf{41}
  (2008), 495--509.

\bibitem{lupton-s-2007b}
G.~Lupton and {S.B.} Smith, \emph{Rationalized evaluation subgroups of a map
  {I}: {S}ullivan models, derivations and {G}-sequences}, J. Pure Appl. Algebra
  \textbf{209} (2007), no.~1, 159--171.

\bibitem{sul-2004}
D.~Sulliivan, \emph{Open and closed string field theory interpreted in
  classical algebaric topology}, Topology, {G}eometry and {Q}uantum {F}ield
  Theory, London Math.Soc. Lecture Notes, vol. 308, Cambridge University Press,
  2004, pp.~344--357.

\bibitem{sul77}
D.~Sullivan, \emph{Infinitesimal computations in topology}, Publ. I.H.E.S.
  \textbf{47} (1977), 269--331.

\bibitem{s-v}
D.~Sullivan and M.~Vigu{\'e}-Poirrier, \emph{The homology theory of the closed
  geodesic problem}, J. Differential Geom. \textbf{11} (1976), 633--644.

\end{thebibliography}
\end{document}